# The L-Functions Have Only Simple Non Trivial Zeros

Tuan Cao-Huu[(1)] and Florin Alan Muscutar[(2)]

[(1)] York University, Glendon College, Toronto, Canada
[(2)] Lorain CCC, 1005 Abbe Road, Eliria, OH44035, USA

*Abstract.* By using an analogy with the case of very close zeros symmetric with respect to the critical line of the Davenport and Heilbronn function, we study the conformal mapping of L-functions in a neighborhood of a hypothetical double zero and conclude that such a zero cannot exist.

**1**. **Introduction**

A zoom on a couple of symmetric zeros with respect to the critical line of the list in [3] (obtained for $t = 520.9438$ with the abscissa $0.51591$ and respectively $0.48409$) gives us an idea of the configuration we should obtain in the case of a double zero of Davenport and Heilbronn function. Indeed, those zeros are enough close to each other and, as a consequence, the configuration of pre-images of rays and circles (see [4]) is similar to that we should obtain for a double zero. The first remark to make is that double zeros for this function, if they exist, should necessarily be located on the critical line. The second is that such a zero is always on the common boundary of two fundamental domains and finally, that one of those domains is bounded by a curve mapped by the function 2:1 onto the negative real half axis, i.e. the respective half axis is a slit for the image of each one of the two fundamental domains. One of those domains is bounded to the right. It can be easily seen that such a configuration should be proper also to double non trivial zeros of any function satisfying the Riemann Hypothesis. Hence, we can use the analogy with the configuration offered by the Davenport and Heilbronn function in order to draw conclusions for this last case. In Fig. 1 below we imagine the two components of the pre-image of the real axis as intersecting each other at the double zero. In this way the component of the pre-image of the negative real half axis becomes the boundary of the fundamental domain from the left, which is hence mapped conformally by the function onto the complex plane with a slit alongside the negative real half axis. When mapping the second fundamental domain we add to that slit the interval $(1,+\infty)$ and possible one or two other arcs starting from $z = 1$.

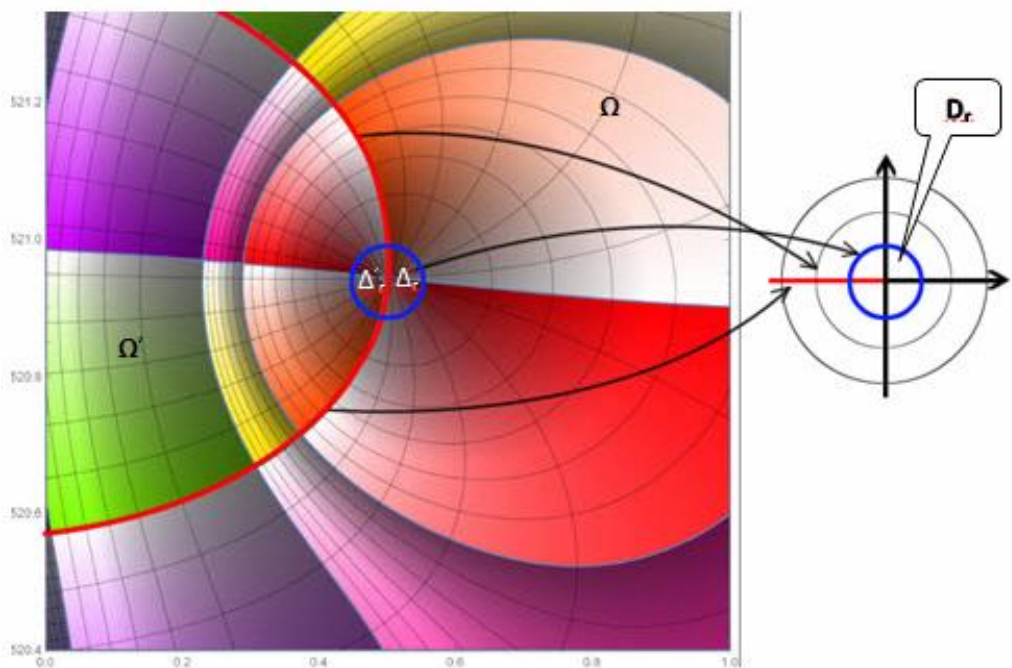

Fig. 1 The local mapping in the neighborhood of a double zero of an L-function

## 2. Functions satisfying the Riemann Hypothesis

The paper [7] deals with the class of functions obtained by analytic continuation of general Dirichlet series of the form:

$$\zeta_{A,\Lambda}(s) = \sum_{n=1}^{\infty} a_n e^{-\lambda_n s} \qquad (1)$$

where $A = \{a_n\}_{n \in \mathbb{N}}$ is a sequence of complex numbers such that $n \to a_n$ is a totally multiplicative function and $\Lambda = \{0 = \lambda_1 \leq \lambda_2 \leq \ldots\}$ is such that $\lim_{n \to \infty} \lambda_n = +\infty$ and the function $n \to \lambda_n$ has the property that for the prime decomposition:

$$n = p_1^{\alpha_1} p_2^{\alpha_2} \ldots p_k^{\alpha_k} \text{ we have } \lambda_n = \alpha_1 \lambda_{p_1} + \alpha_2 \lambda_{p_2} + \ldots + \alpha_k \lambda_{p_k} \qquad (2)$$

Let us denote by $\sigma_c$ the abscissa of convergence of the series (1). Then, as shown in [7], if $n \to a_n$ is a totally multiplicative function and $n \to \lambda_n$ satisfies (2), then

$$\zeta_{A,\Lambda}(s) = \Pi_p (1 - a_p e^{-\lambda_p s}) \qquad (3),$$

where $p$ runs through prime numbers and the series (3) converges uniformly on compact sets in the half plane $\operatorname{Re} s > \sigma_c$. Also, if $\sigma_c < 1/2$ and the function $\zeta_{A,\Lambda}(s)$ satisfies a Riemann-type of functional equation, then the non trivial zeros of $\zeta_{A,\Lambda}(s)$ are all on the critical line. We will use next the Fig. 1 as a guide in order to describe the mapping in the neighborhood of a hypothetical double non trivial zero of $\zeta_{A,\Lambda}(s)$. For simplicity, we denote by $f(s)$ such a function.

In what follows, we are using techniques perfected by D. Ghisa in his monograph [4] and a series of papers (see [2], [5], [6], [7]). For any function obtained by analytic continuation of a series (1) the complex plane can be partitioned into infinitely many strip $S_k$, $k \in \mathbb{Z}$, $0 \in S_0$ and $S_k$ is below $S_{k+1}$ for every $k$ (*Ghisa's strips*) such that every strip $S_k$ is mapped by the function (not necessarily one to one) onto the whole complex plane with a slit alongside the interval $(1, +\infty)$ of the real axis. The boundary of every $S_k$ is formed by two unbounded curves $\Gamma'_k$ and $\Gamma'_{k+1}$, the counting being in ascending order. Every strip $S_k$, $k \neq 0$ contains a unique unbounded component of the pre-image by $\zeta_{A,\Lambda}(s)$ of the unit disc, as well as a unique component $\Gamma_{k,0}$ of the pre-image by $\zeta_{A,\Lambda}(s)$ of the real axis, which is mapped by $\zeta_{A,\Lambda}(s)$ one to one onto the interval $(-\infty, 1)$ of the real axis. If we denote by $j_k$ the number of zeros of $\zeta_{A,\Lambda}(s)$ in $S_k$, $k \neq 0$, then there are $j_k - 1$ other components $\Gamma_{k,j}$, $j \neq 0$ of the pre-image of the real axis in $S_k$ which are mapped one to one onto the whole real axis. The curves $\Gamma'_k$ and $\Gamma_{k,j}$ have been known for a long time, but only recently D. Ghisa revealed properties of these curves which make them a valuable tool in the study of the functions (1).

The derivative $\zeta'_{A,\Lambda}(s)$ has properties similar to those of $\zeta_{A,\Lambda}(s)$, except that $\lim_{\sigma \to +\infty} \zeta_{A,\Lambda}(\sigma + it) = 1$, while $\lim_{\sigma \to +\infty} \zeta'_{A,\Lambda}(s) = 0$, which makes that the corresponding curves $\Upsilon_{k,0}$ contain no zero of $\zeta'_{A,\Lambda}(s)$. For every $k \neq 0$ and for every $j$ the curves $\Upsilon_{k,j}$ and $\Upsilon'_k$ are such that every $\Upsilon_{k,j}$ intersects a unique curve $\Gamma_{k,j}$, respectively every $\Upsilon'_k$ intersects a unique curve $\Gamma'_k$ at points where the tangent to the last ones is horizontal. These are *Ghisa's intertwined curves* and the intertwining of these curves happens in a very precise way, namely respecting *Ghisa's color matching rule* (see [5] and [6]).

Both $\zeta_{A,\Lambda}(s)$ and $\zeta'_{A,\Lambda}(s)$ obey yet another law, namely the *color alternating rule.* These rules allow a partition of every strip $S_k$ into subsets such that the interior $\Omega_{k,j}$ of each one of them is a fundamental domain in the sense defined by Ahlfors (see in [1], page 98). We will call $\Omega_{k,j}$ *Ghisa's fundamental domains*. They play a crucial role in the study of the distribution of zeros of these functions, as well as in answering questions related to the multiplicity of those zeros.

For simplicity, we will use the generic notation $f(s)$ for a function (1). Suppose that $s_0$ is a double zero of $f(s)$, i.e. $f(s_0) = 0$ and $f'(s_0) = 0$, $f''(s_0) \neq 0$, i.e. $s_0$ is a simple zero for $f'(s)$. Let $\Omega$ and $\Omega'$ be the fundamental domains having $s_0$ as common boundary point. Then there are two curves $\Gamma_{k,j}$ and $\Gamma_{k,j'}$ passing through $s_0$. Let $\Upsilon_{k,j}$ and $\Upsilon_{k,j'}$ be the intertwined curves (see [6]). If $j \neq 0$ and $j' \neq 0$, the curves $\Upsilon_{k,j}$ and $\Upsilon_{k,j'}$ must intersect each other, since their intertwined curves, none of which is $\Gamma_{k,0}$ intersect. But this can happen only at $s_0$ (see [8]), which would imply that $s_0$ is a double zero for $f'(s)$, and this is a contradiction, therefore $j = 0$ or $j' = 0$. Suppose that $j' = 0$.

Now $\Upsilon_{k,j'}$ does not pass through $s_0$. Let $L$ and $L'$ be the slits corresponding to $\Omega$ and $\Omega'$. Having in view the analogy with the configuration corresponding to a double zero for the Davenport and Heilbronn function, we know that $L'$ is the negative real half axis. The boundary $\partial \Omega'$ is mapped by $f(s)$ (two to one) onto $L'$. Let us denote $\Omega_\varphi = \Omega \cup \Omega' \cup \partial \Omega'$ and notice that $\Omega_\varphi$ is a simply connected domain. Let us define $\varphi : \Omega \to \Omega'$ by

$$\varphi(s) = f_{|\Omega'}^{-1}(f_{|\Omega}(s)) \qquad (4)$$

and notice that $f_{|\Omega}$ and $f_{|\Omega'}$ extend naturally to $\Omega \cup (\partial \Omega' \cap \Gamma_{k,0})$, respectively $\Omega' \cup (\partial \Omega' \cap \Gamma_{k,-1})$ and if $s \in \partial \Omega' \cap \Gamma_{k,0}$, or $s \in \partial \Omega' \cap \Gamma_{k,-1}$ then $f_{|\Omega}(s)$ and respectively $f_{|\Omega'}(s)$ should be seen as symmetric points on different borders of $L'$. With this convention we can extend $\varphi(s)$ to an involution of $\Omega_\varphi$ by taking $\varphi(\varphi(s)) = s$ for every $s \in \Omega_\varphi$. The function $\varphi(s)$ is a conformal mapping of $\Omega_\varphi$ onto itself having the fixed point $s_0$. Moreover, the equality $\varphi'(\varphi(s))\varphi'(s) = 1$ implies that $\varphi'^2(s_0) = 1$. Let $g : \mathbb{C} \setminus L \to \Omega$ be the inverse function of $f_{|\Omega \cup (\partial \Omega' \cap \Gamma_{k,0})}(s)$. Then, for every $z \in \mathbb{C} \setminus L$ we have

$$f(\varphi(g(z))) = z \qquad (5)$$

The functions in (5) are all analytic in their domains and differentiating with respect to $z$ we get

$$f'(\varphi(g(z)))\varphi'(g(z))g'(z) = 1 \qquad (6)$$

This equality holds in a neighborhood of $0$ except at the points of the negative real half axis, where $g'(z)$ is not defined. However $g(z)$ is defined and $g(0) = s_0$.

We keep in mind that $\varphi(s_0) = s_0$ and $f'(\varphi(g(0))) = f'(s_0) = 0$. With the notation $s = g(z)$ and having in view the expression of the derivative of the inverse function, the equation (6) can be written as:

$$f'(\varphi(s))\varphi'(s) = f'(s) \qquad (7)$$

Since $f'(s_0) = f'(\varphi(s_0)) = 0$ and $|\varphi'(s)| = 1$, the equation (7) expresses the fact that

$|f'(s)|$ decreases to zero with the same rate in $\Omega$ and $\Omega'$ as $s \to s_0$.

If $D_r$ is a disc centered at the origin of radius $r$ and $\Delta_r$, and $\Delta'_r$ are the components of its pre-image included in $\Omega$, respectively $\Omega'$, then

$$\int_{\Delta_r} |f'(s)|^2 d\sigma dt = \int_{\Delta'_r} |f'(s)|^2 d\sigma dt = \pi r^2, \qquad (8)$$

both integrals representing the area of the image by $f(s)$ of $\Delta_r$, respectively $\Delta'_r$. One of the integrals in (8) can be obtained from the other also by the change of variable defined by (7).

Moreover, if we denote by $\gamma$ and $\gamma'$ the component included in $\Omega$, respectively $\Omega'$ of the pre-image of the circle centred at the origin and of radius $r$, then:

$$\int_{\gamma}|f'(s)||ds| = \int_{\gamma'}|f'(s)||ds| = 2\pi r \qquad (9)$$

both integrals representing the length of the image by $f(s)$ of the curves $\gamma$, respectively $\gamma'$. Also, the integrals of $|f'(s)|$ on pre-image of rays of such a circle are all equal, since each one is the length of such a ray.

At the infinitesimal scale, $\Delta_r$ and $\Delta'_r$ can be considered as half discs of the radius $\sqrt{r}$ and $f(s) = (s - s_0)^2 h(s)$, where $h(s_0) \neq 0$ (see [1], page 133). Indeed, for such a function the equations (8) and (9) are satisfied and we also have $\varphi(s_0 + \sigma + it) \approx s_0 - \sigma - it$ and $\varphi$ is an involution of $\Delta_r \cup \Delta'_r$ with the fixed point $s_0$.

**Theorem.** *Suppose that $n \to a_n$ is a totally multiplicative function and $n \to \lambda_n$ satisfies the equation* (2). *If $\sigma_c < 1/2$ and $\zeta_{A,\Lambda}(s)$ satisfies a Riemann-type of functional equation, then all the non trivial zeros of $\zeta_{A,\Lambda}(s)$ are simple.*

*Proof:* Let us denote as in [7] $B \subset \Omega$ and $B' \subset \Omega'$ the components of the pre-image by $\zeta_{A,\Lambda}(s)$ of $\mathbb{C} \setminus (L \cup L')$ and $\Omega_c = (B \cup B') \cap [\{s \mid \operatorname{Re} s > \sigma_c\} \cup \{s \mid \operatorname{Re} \varphi(s) > \sigma_c\}]$ and define $f_n(s)$ in $\Omega_c$ by

$$f_n(s) = \Pi_{p \leq n}(1 - a_p e^{-\lambda_p s}), \text{ in } (B \cup B') \cap \{s \mid \operatorname{Re} s > \sigma_c\} \qquad (10)$$

$$f_n(s) = \Pi_{p \leq n}(1 - e^{-\lambda_p \varphi(s)}), \text{ in } (B \cup B') \cap \{s \mid \operatorname{Re} s \leq \sigma_c\} \cap \{s \mid \operatorname{Re} \varphi(s) > \sigma_c\}$$

where $p$ runs through prime numbers. Then $f_n(s)$ are analytic in $\Omega_c$ which includes a neighborhood of $s_0$. Moreover, $\lim_{n \to \infty}[1 / f_n(s)] = \zeta_{A,\Lambda}(s)$ uniformly on compact sets included in $\Omega_c$. As in [7], we check easily that the equation $1 - a_p e^{-\lambda_p s} = 0$ cannot have any solution in $\Omega_c$. This fact guarantees that $\lim_{n \to \infty} f_n(\varphi(s)) / f_n(s) = \zeta_{A,\Lambda}(s) / \zeta_{A,\Lambda}(\varphi(s)) = 1$ in $\Omega_c$.

Let us notice that if $s = s_0 + \sigma + it$, we can take $\sigma$ such that $\operatorname{Re}(s - \varphi(s)) > 0$. Then

$$f_n(\varphi(s)) / f_n(s) = \Pi_{p \leq n}[(e^{\lambda_p \varphi(s)} - a_p) / (e^{\lambda_p s} - a_p)]e^{\lambda_p(s - \varphi(s))} \qquad (11)$$

and the sequence (11) is divergent since $e^{\lambda_p(s - \varphi(s))} \to \infty$ as $p \to \infty$ and none of the factors $e^{\lambda_p \varphi(s)} - a_p$ can be zero in order to offset that limit. This contradicts the convergence of the sequence (10) at $s$ or at $\varphi(s)$. The final conclusion is that $\zeta_{A,\Lambda}(s)$ cannot have any double zero.¤

This analysis has been performed for a hypothetical double zero of $\zeta_{A,\Lambda}(s)$ satisfying the hypothesis of the theorem.

What about the Davenport and Heilbronn function?

We have already noticed that if a double zero existed for this function then its location must have been on the critical axis. Yet, for such a zero, the same contradiction can be drawn as in the previous theorem. Therefore, the Davenport and Heilbronn function has no double non trivial zero.

A similar analysis can be performed by using [1], page 133 for a zero of a higher order of multiplicity.